\documentclass[a4paper,11pt,mathcal]{article}
\usepackage{amsfonts,eucal,amssymb,amsmath,amsthm,amscd} 

\newtheorem{theorem}{Theorem}[section]
\newtheorem{corollary}[theorem]{Corollary}

\newtheorem{proposition}[theorem]{Proposition}
\newtheorem{lemma}[theorem]{Lemma}
\newtheorem{definition}[theorem]{Definition}
\newtheorem{question}[theorem]{Question}

\newenvironment{Pf}{\medskip \noindent {\bf Proof: }}
   {$\diamondsuit$ }

\newcommand{\Hom}{\text{Hom}}
\newcommand{\Mon}{\text{Mon}}
\newcommand{\Ad}{\text{Ad }}
\newcommand{\AdgP}{\text{Ad}_{\mathfrak{g}} P}
\newcommand{\Aut}{\text{Aut }}

\newcommand{\Kill}{\text{Kill}}

\newcommand{\Gr}{\text{Gr }}
\newcommand{\Zar}{\text{Zar}}

\newcommand{\SL}{\text{SL}}
\newcommand{\GL}{\text{GL}}

\newcommand{\liep}{\ensuremath{\mathfrak{p}}}

\newcommand{\lieg}{\ensuremath{\mathfrak{g}}}
\newcommand{\lien}{\ensuremath{\mathfrak{n}}}

\newcommand{\lieh}{\ensuremath{\mathfrak{h}}}
\newcommand{\lies}{\ensuremath{\mathfrak{s}}}

\newcommand{\lieu}{\ensuremath{\mathfrak{u}}}
\newcommand{\liec}{\ensuremath{\mathfrak{c}}}

\newcommand{\RR}{\ensuremath{{\bf R}}}

\newcommand{\NN}{\ensuremath{{\bf N}}}

\newcommand{\dduzero}{\ensuremath{\left.\frac{\mbox{d}}{\mbox{d}u} \right|_0}}
\newcommand{\ddszero}{\ensuremath{\left.\frac{\mbox{d}}{\mbox{d}s} \right|_0}}
\newcommand{\ddtzero}{\ensuremath{\left.\frac{\mbox{d}}{\mbox{d}t} \right|_0}}
\newcommand{\kddtzero}{\left.\frac{\mbox{d}^k}{\mbox{d}t^k} \right|_0}
\newcommand{\nddszero}{\left.\frac{\mbox{d}^n}{\mbox{d}s^n} \right|_0}

\begin{document}
\pagenumbering{arabic}

\title{A Frobenius theorem for Cartan geometries, with applications}
\author{Karin Melnick\footnote{karin.melnick@yale.edu,
    partially supported by NSF fellowship DMS-855735}}
\date{\today}
\maketitle

\section{Introduction}

The classical result on local orbits in geometric manifolds is
Singer's homogeneity theorem for Riemannian manifolds \cite{singer.locisom}: given a Riemannian
manifold $M$, there exists $k$, depending on $\dim M$, such that if every $x,y \in M$
are related by an infinitesimal isometry of order $k$, then $M$ is locally homogeneous.
An open subset $U \subseteq M$ of a geometric manifold is \emph{locally
  homogeneous} if for every $x,x' \in U$, there is a \emph{local automorphism}
$f$ in $U$ with $f(x) = x'$.  Such a local automorphism is a diffeomorphism from a
neighborhood $V$ of $x$ in $U$ to a neighborhood of $x'$ in $U$, with $f$
an isomorphism between the geometric structures restricted to $V$ and $f(V)$.

Gromov extended Singer's theorem to manifolds with \emph{rigid geometric structures
  of algebraic type} in \cite[1.6.G]{gromov.rgs}.  He also proved the celebrated
open-dense theorem (\cite[3.3.A]{gromov.rgs}) and a stratification for orbits
of local automorphisms of
such structures on compact real-analytic manifolds (see \cite[3.4]{gromov.rgs} and
\cite[3.2.A]{dag.rgs}).  The open-dense theorem says that if $M$ is a smooth
manifold with smooth rigid geometric structure of algebraic type, and if there
is an orbit for local automorphisms that is dense in $M$, then $M$ contains an open, dense, locally homogeneous subset.
A crucial ingredient for Gromov's theorems is his difficult Frobenius
theorem, which says that infinitesimal isometries of sufficiently high order
can be integrated to local isometries near any point on a real-analytic
manifold, and near \emph{regular} points in the smooth case.

This article treats \emph{Cartan geometries}, a notion of geometric structure
less flexible than Gromov's rigid geometric structures, but
still including essentially all classical geometric structures with
finite-dimensional automorphism groups, such as
pseudo-Riemannian metrics, conformal pseudo-Riemannian structures in dimension at least $3$, and a broad class of CR
structures.  The central result is a Frobenius theorem for Cartan geometries
(\ref{thm.can.frobenius}, \ref{thm.smooth.regular}),
which is considerably easier in this setting, and is in fact broadly modeled
on the paper \cite{nomizu.killing.gens} of Nomizu from 1960 treating Riemannian isometries (see also \cite{cqb.centralizer}).

From the Frobenius theorem we obtain the stratification and open-dense
theorems as in \cite{gromov.rgs} for local Killing fields of Cartan geometries (\ref{thm.can.stratification}, \ref{thm.smooth.opendense}).  The embedding theorem for automorphism groups of Cartan geometries proved in \cite{bfm.zimemb}, combined with the Frobenius theorem, gives rise to centralizer and Gromov representation theorems for real-analytic Cartan geometries (\ref{thm.centralizer}, \ref{thm.gromovrep}), which can be formulated for actions that do not necessarily preserve a finite volume.





A Cartan geometry infinitesimally models a manifold on a homogeneous space.

\begin{definition}
\label{def.cg}
A \emph{Cartan geometry} on a manifold $M$ modeled on a homogeneous space $G/P$ is a triple $(M,B,\omega)$ where 
$B$ is a principal $P$-bundle over $M$, and $\omega$ is a $\lieg$-valued $1$-form on $B$ satisfying
\begin{enumerate}
\item{$\omega_b : T_bB \rightarrow \lieg$ is a linear isomorphism for all $b \in B$}
\item{for all $X \in \liep$, if $X^\ddag$ is the fundamental vector field on $B$ corresponding to $X$, then $\omega_b(X^\ddag) = X$ at all $b \in B$.}
\item{$R_g^* \omega = \Ad g^{-1} \circ \omega$ for all $g \in P$}
\end{enumerate}
\end{definition}

\begin{definition}
Let $(M,B,\omega)$ be a Cartan geometry.  An \emph{automorphism} of $(M,B,\omega)$ is a diffeomorphism $f$ of $M$ that lifts to a bundle automorphism $\tilde{f}$ of $B$ satisfying $\tilde{f}^* \omega = \tilde{f}$.
\end{definition}

Let $(M,B,\omega)$ be a Cartan geometry modeled on $G/P$.  We will make the following standard assumptions on $G/P$:
\begin{enumerate}
\item{$G$ is connected.} 
\item{$P$ contains no nontrivial normal subgroup of $G$.  (Suppose that $N
    \lhd G$ were such a subgroup.  Then let $G' = G/N$ and $P' = P/N$.  If $(M,B,\omega)$ is a Cartan geometry modeled on $(\lieg,P)$, then $\omega$ descends to a $\lieg'$-valued $1$-form on $B'= B/N$, giving a Cartan geometry $(M,B/N,\omega')$ modeled on $G'/P'$.)}
\item{$P$ is an analytic subgroup of $G$.}
\end{enumerate}

In section \ref{section.gromovrep} we will further assume that $\AdgP$ is an algebraic subgroup of $\Aut \lieg$.  In this case, the Cartan geometry $(M,B,\omega)$ is said to be \emph{algebraic type}.

\emph{Acknowledgements:}  I thank Charles Frances, Gregory Margulis, Amir
Mohammadi, and especially David Fisher, for
helpful conversations during the writing of this paper.

\section{Baker-Campbell-Hausdorff formula}

The main proposition of this section asserts that the usual BCH formula
holds to any finite order with $\omega$-constant vector fields
on $B$ in place of left-invariant vector fields.  When $(M,B,\omega)$ is real-analytic, this formula gives the Taylor
series at each point of $b$ for the flow along two successive $\omega$-constant vector
fields, in terms of the exponential coordinates.

For $X,Y \in \lieg$, define 
\begin{eqnarray*}
\alpha & : & \lieg \times \lieg \rightarrow \lieg \\
\alpha & : & (X,Y) \mapsto \log_e ( \exp_e X \cdot \exp_e Y)
\end{eqnarray*}

where $\exp_e$ is the group exponential map $\lieg \cong T_e G \rightarrow G$,
and $\log_e$ the inverse of $\exp_e$.  The exponential map of $G$ can be
considered a function $G \times \lieg \rightarrow G$, with 
$$\exp(g,X) = \exp_g X = g \cdot \exp X$$

It is the flow for time $1$ with initial value $g$ along the left-invariant vector field
corresponding to $X$.  Note that 
$$ \exp X \cdot \exp Y = \exp(\exp(e,X),Y)$$

For any $k \in \NN$, there exist functions
$a_1, \ldots, a_k,$ and $R$ of $(X,Y)$ such that
$$ \alpha(tX,tY) = ta_1(X,Y) + \cdots + \frac{t^k}{k!}a_k(X,Y) + t^k R(tX,tY)$$
where 
$$\lim_{t \rightarrow 0} R(tX,tY) = 0$$

These functions are given by the BCH formula, and they are rational multiples
of iterated brackets
of $X$ and $Y$.  For example, 
\begin{eqnarray*}
 a_1(X,Y) & = & X + Y \\
a_2(X,Y) & = & [X,Y]
\end{eqnarray*}
and 
$$ a_3(X,Y) = \frac{1}{2}([X,[X,Y]]+[Y,[Y,X]])$$

For any Lie algebra $\lieu$, not necessarily finite-dimensional, with a linear
injection $\rho: \lieg \rightarrow \lieu$, the functions $a_k$ define obvious
functions $ a_k : \rho(\lieg) \rightarrow \lieu$, evaluated by taking iterated brackets
in $\lieu$.

In the bundle $B$ of the Cartan geometry, denote by $\exp$ the
exponential map $B \times \lieg \rightarrow B$, defined on a neighborhood of
$B \times \{ 0 \}$ and by $\log_b$ the inverse of $\exp_b$, defined on a
normal neighborhood of $b$.
For any $b \in B$, define, for sufficiently small $X,Y \in \lieg$
$$\zeta_b(X,Y) = \log_b(\exp(\exp(b, X),Y))$$

As above, there exist functions $z_1, \ldots, z_k$, corresponding to the time derivatives of
$\zeta_b(tX,tY)$ up to order $k$, and a remainder function.

\begin{proposition}
\label{prop.bch}
Let $G$ be a Lie group with Lie algebra $\lieg$ and $(M,B,\omega)$ a Cartan geometry modeled on a
homogeneous space of $G$.  Let $a_k$ and $z_k$ be the coefficients of $t^k/k!$ in
the respective order-$k$ Taylor approximations of the above functions $\alpha$
and $\zeta_b$.  Then
$$ z_k(X,Y) = \omega_b(a_k(\widetilde{X},\widetilde{Y}))$$

where $\widetilde{X}$ and $\widetilde{Y}$ are the $\omega$-constant vector fields on $B$
corresponding to $X$ and $Y$, respectively.
\end{proposition}

\begin{Pf}
Fix $X,Y \in \lieg$, and let $Z(t) = \zeta_b(tX,tY)$.  The following lemmas
give two different ways to compute, for an arbitrary $C^k$ function $\varphi$
on $B$ and $b \in B$, the derivative 
$$ \kddtzero \varphi(\exp(b, Z(t)))$$

\begin{lemma}
\label{kthder1.lemma}
For $X \in \lieg$, $b \in B$, and $\varphi \in C^k(B)$, 
$$ \kddtzero \varphi(\exp(b,tX)) = \left. \widetilde{X}^k.\varphi\right|_b$$
\end{lemma}

\begin{Pf}
For $k = 1$, 
$$ \ddtzero \varphi(\exp(b,tX)) = \varphi_{*b}((\exp_b)_*(X)) = \left.\widetilde{X}.\varphi\right|_b$$

Now let $n \geq 1$ and suppose that the formula holds for all $k \leq n$.  Then
\begin{eqnarray*}
\left.\widetilde{X}^{n+1}.\varphi \right|_b & = & \left.\widetilde{X}.\widetilde{X}^n.\varphi\right|_b \\
 & = & \ddtzero(\widetilde{X}^n.\varphi)(\exp(b,tX)) \\
& = & \ddtzero \nddszero \varphi(\exp(\exp(b,tX),sX)) \\
& = & \ddtzero \nddszero \varphi(\exp(b,(t+s)X)) \\
& = & \left. \frac{\mbox{d}^{n+1}}{\mbox{d}u^{n+1}} \right|_0 \varphi(\exp(b,uX))
\end{eqnarray*}

where $u = t+s$. \end{Pf}

\begin{corollary}
\label{kthder1.cor}
For $X,Y \in \lieg$,
$$ \kddtzero \varphi(\exp(\exp(b,tX),tY)) = \sum_{m+n=k} \frac{k!}{m!n!}
\left.\widetilde{X}^n. \widetilde{Y}^m.\varphi \right|_b$$
\end{corollary}

\begin{Pf}
By two applications of lemma \ref{kthder1.lemma},
\begin{eqnarray*}
\left.\widetilde{X}^n.\widetilde{Y}^m.\varphi \right|_b & = & \nddszero (\widetilde{Y}^m.\varphi)(\exp(b,sX)) \\
& = & \nddszero \left. \frac{\mbox{d}^m}{\mbox{d}t^m} \right|_0 \varphi(\exp(\exp(b,sX),tY))
\end{eqnarray*}
The desired formula follows.
\end{Pf}

\begin{lemma}
\label{kthder2.lemma}
Let $Z(t)$ be a curve in $\lieg$, $b \in B$, and $\varphi \in C^k(B)$.  Then
$$ \kddtzero \varphi(\exp_b Z(t)) = \kddtzero \left[ \sum_{n=0}^k \frac{1}{n!}
\left.[t\widetilde{Z'(0)} + \cdots + \frac{t^k}{k!}
  \widetilde{Z^{(k)}(0)}]^n.\varphi\right|_{\exp_b Z(0)} \right] $$
where $\widetilde{Z^{(l)}(0)}$ is the $\omega$-constant vector field on $B$
  evaluating to
$$ \omega((\exp_b)_{* Z(0)} Z^{(l)}(0))$$ 
\end{lemma}

\begin{Pf}
Let $c(t) = \exp_b Z(t)$. Equality is clear when $k=0$.  When $k=1$, the left side is
$$\left.\widetilde{Z'(0)}.\varphi\right|_{c(0)}$$
and the right side is
$$ \ddtzero \left[ \varphi(c(0)) + \left.t\widetilde{Z'(0)}.\varphi
  \right|_{c(0)} \right] = \left. \widetilde{Z'(0)}.\varphi \right|_{c(0)}$$

Now let $n \geq 1$, and suppose that the formula holds for $k \leq n$ for any
curve $Z(t)$.  Then

\begin{eqnarray*}
 \left. \frac{\mbox{d}^{n+1}}{\mbox{d}t^{n+1}}\right|_0 \varphi(c(t))  & = & \ddtzero \left. \frac{\mbox{d}^n}{\mbox{d}s^n} \right|_t \varphi(c(s)) \\
& = & \ddtzero \nddszero \varphi(c(s+t)) \\
& = & \ddtzero \nddszero \left[ \sum_{k=0}^n \frac{1}{k!} \left.[s\widetilde{Z'(t)} +
  \cdots + \frac{s^n}{n!} \widetilde{Z^{(n)}(t)}]^k.\varphi \right|_{c(t)} \right] 
\end{eqnarray*}

where $\widetilde{Z^{(l)}(t)}$ is the $\omega$-constant vector
field evaluating to 
$$ \omega_b((\exp_b)_{*Z(t)} Z^{(l)}(t))$$
at $c(t)$.  Continuing, the last expression equals

\begin{eqnarray*}
&   & \nddszero \ddtzero \left[ \sum_{k=0}^{n+1} \frac{1}{k!} \left.[s\widetilde{Z'(t)} +
  \cdots + \frac{s^{n+1}}{(n+1)!}  \widetilde{Z^{(n+1)}(t)}]^k.\varphi\right|_{c(t)} \right] \\
& = & \nddszero \left[ \sum_{k=1}^{n+1} \frac{1}{(k-1)!} [s\widetilde{Z'(0)} +
  \cdots + \frac{s^{n+1}}{(n+1)!}
  \widetilde{Z^{(n+1)}(0)}]^{k-1}.[s\widetilde{Z''(0)} \right. \\
 & \qquad & +  \cdots + \left. \frac{s^n}{n!} \widetilde{Z^{(n+1)}(0)} + \frac{s^{n+1}}{(n+1)!} \widetilde{Z^{(n+2)}(0)}].\varphi\right|_{c(0)} \\
& \qquad & + \left. \sum_{k=0}^{n+1} \frac{1}{k!} \left. \widetilde{Z'(0)}.[s\widetilde{Z'(0)} +
  \cdots + \frac{s^{n+1}}{(n+1)!} \widetilde{Z^{(n+1)}(0)}]^k.\varphi\right|_{c(0)} \right] \\
\end{eqnarray*}

using that 
\begin{eqnarray*}
\ddtzero \left. \widetilde{Z^{(m_1)}(t)}. \ldots
  .\widetilde{Z^{(m_l)}(t)}.\varphi \right|_{c(t)} & = & \sum_{i=1}^l
\left. \widetilde{Z^{(m_1)}(0)}. \ldots .\widetilde{Z^{(m_i+1)}(0)}. \ldots. \widetilde{Z^{(m_l)}(0)}.\varphi
\right|_{c(0)} \\
&  & + \left. \widetilde{c'(0)}.\widetilde{Z^{(m_1)}(0)}. \ldots
  .\widetilde{Z^{(m_l)}(0)}.\varphi \right|_{c(0)}
\end{eqnarray*}

Continuing, we have
\begin{eqnarray*}
& = & \nddszero \left[ \sum_{k=1}^{n+1} \frac{1}{(k-1)!} [s\widetilde{Z'(0)} +
  \cdots + \frac{s^{n+1}}{(n+1)!}  \widetilde{Z^{(n+1)}(0)}]^{k-1}.[\widetilde{Z'(0)} + s\widetilde{Z''(0)} \right. \\
& \qquad & \left. \left. + \cdots + \frac{s^n}{n!} \widetilde{Z^{(n+1)}(0)}].\varphi\right|_{c(0)} \right] \\
& = & \nddszero \frac{\mbox{d}}{\mbox{d}s} \left[ \sum_{k=0}^{n+1}
\left. \frac{1}{k!} [s\widetilde{Z'(0)} + \cdots + \frac{s^{n+1}}{(n+1)!}
  \widetilde{Z^{(n+1)}(0)}]^k.\varphi\right|_{c(0)} \right] \\
& = & \left. \frac{\mbox{d}^{n+1}}{\mbox{d}s^{n+1}} \right|_0 \left[
  \sum_{k=0}^{n+1} \left. \frac{1}{k!} [s\widetilde{Z'(0)} + \cdots +
  \frac{s^{n+1}}{(n+1)!} \widetilde{Z^{(n+1)}(0)}]^k.\varphi\right|_{c(0)} \right]
\end{eqnarray*}
\end{Pf}

Now 
$$ \exp(\exp(b,tX),tY) = \exp(b,\zeta_b(tX,tY)) = \exp_b(Z(t))$$

for $Z(t) = \zeta_b(tX,tY)$.  Note that $Z^{(k)}(0) = z_k(X,Y)$.  Corollary \ref{kthder1.cor} gives
$$ \kddtzero \varphi(\exp_b Z(t)) = \sum_{m+n=k} \frac{k!}{m!n!}
\left.\widetilde{X}^n. \widetilde{Y}^m.\varphi \right|_b$$

On the other hand, lemma \ref{kthder2.lemma} gives for the same derivative
$$ \kddtzero \varphi(\exp_b Z(t)) = \kddtzero \left[ \sum_{n=0}^k \frac{1}{n!}
\left.[t\widetilde{Z'(0)} + \cdots + \frac{t^k}{k!}
  \widetilde{Z^{(k)}(0)}]^n.\varphi\right|_b \right] $$

With these two formulas, the coefficients $Z^{(k)}(0) = z_k(X,Y)$ can be
recursively computed in terms of products of $\widetilde{X}$ and
$\widetilde{Y}$.  Of course, these formulas hold in the group $G$ with the
usual exponential map, so
they yield the same expressions, actually involving brackets of $\widetilde{X}$
and $\widetilde{Y}$, for $a_k(X,Y)$ and $z_k(X,Y)$.
\end{Pf}

\section{Frobenius theorem}
\label{section.frobenius}

Throughout this section, $(M,B,\omega)$ is a Cartan geometry modeled on
$G/P$.  Soon we will impose the assumption that $(M,B,\omega)$ is $C^\omega$.  The curvature of a Cartan geometry is a $\lieg$-valued $2$-form on
$B$ defined by
$$ \Omega(X,Y) = \mbox{d}\omega(X,Y) + [\omega(X),\omega(Y)]$$
If $X \in \liep$, then $\Omega(X,Y)$ vanishes \cite[5.3.10]{sharpe}.  Let 
$$ V = (\wedge^2 (\lieg / \liep)^*) \otimes \lieg$$
The form $\omega$ gives an identification $TB \cong B \times \lieg$, under
which the curvature corresponds to
a function $ K : B \rightarrow V$
$$ K : b \mapsto (\omega_b^{-1} \circ \sigma)^* \Omega_b $$
where $\sigma$ is any linear section $\lieg / \liep \rightarrow \lieg$.

The group $P$ acts on $V$ linearly by 
$$(p. \varphi)(u,v) = (\Ad p \circ \varphi)((\bar{\Ad} p^{-1}) u,(\bar{\Ad} p^{-1})
v)$$
where $\bar{\Ad}$ is the quotient representation of $\Ad P$ on $\lieg / \liep$.
The curvature map is $P$-equivariant \cite[5.3.23]{sharpe}:
$$ K(bp^{-1}) = p.K(b)$$

For $m \in \NN$, define the \emph{$\omega$-derivative} of order $m$ of $K$
\begin{eqnarray*}
D^m K & : & B \rightarrow \Hom(\otimes^m \lieg ,V) \\
D^m K(b) & : & X_1 \otimes \cdots \otimes X_m \mapsto (\widetilde{X}_1 \ldots \widetilde{X}_m. K)(b)
\end{eqnarray*}
where, as above, $\widetilde{X}$ is the $\omega$-constant vector field on $B$
with value $X$.
Note that $D^mK(b)$ is not a symmetric homomorphism, because the
$\omega$-constant vector fields $\widetilde{X}$ do not come from coordinates on
$B$.  Neither can it be interpreted as a tensor on $B$, because $D^m K(b)$ is
not linear over the ring of functions $C^\infty(B)$.  It does suffice,
however, to determine the
$m$-jet $j^m_b K$, because any vector field on $B$ is a $C^\infty(B)$-linear
combination of $\omega$-constant vector fields.

\begin{proposition}
The $\omega$-derivative is $P$-equivariant for each $m \geq 0$:
$$ D^mK(bp^{-1}) = p \circ D^m K(b) \circ \Ad^m p^{-1}$$
where $\Ad^m$ is the tensor representation on $\otimes^m \lieg$ of $\Ad P$.
\end{proposition}

\begin{Pf}
The assertion holds for $m = 0$ by the equivariance of $K$ cited above.
Suppose it holds for all $m \leq r$.  Then for any $X_1, \ldots, X_{r+1} \in \lieg$,
\begin{eqnarray*}
& & (\widetilde{X}_1 \ldots \widetilde{X}_{r+1}. K)(bp^{-1})  =  
\ddtzero (\widetilde{X}_2 \ldots \widetilde{X}_{r+1}. K)(\exp(bp^{-1},tX_1)) \\
& = & \ddtzero (\widetilde{X}_2 \ldots \widetilde{X}_{r+1}. K)(\exp(b,(\Ad p^{-1}) t X_1)
p^{-1}) \\
& = & \ddtzero p.(D^rK(\exp(b,(\Ad p^{-1})tX_1)) ((\Ad p^{-1}) X_2, \ldots, (\Ad
p^{-1}) X_{r+1})) \\
& = & p.(D^{r+1} K(b)((\Ad p^{-1})X_1, \ldots, (\Ad p^{-1}) X_{r+1}))
\end{eqnarray*}
so by induction it is true for all $m \geq 0$.
\end{Pf}

\begin{definition}
For $m \geq 1$, two points $b,b'$ of $B$ are \emph{$m$-related} if 
$$ D^r K(b) = D^r K(b')$$
for all $1 \leq r \leq m$.  They are \emph{$\infty$-related} if they are $m$-related for all $m$.
\end{definition}

For $\varphi \in \Hom(\otimes^r \lieg, V)$ and $X \in \lieg$, the contraction
$\varphi \llcorner X \in \Hom(\otimes^{r-1} \lieg, V)$ is given by 
$$ (\varphi \llcorner X)(X_1, \ldots, X_{r-1}) = \varphi(X, X_1, \ldots, X_{r-1})$$

\begin{definition}
For $m \geq 1$, the \emph{Killing generators of order $m$} at $b \in B$, denoted $\Kill^m(b)$, comprise all $A
\in \lieg$ such that, for all $1 \leq r \leq m$, the contraction
$$ D^r K(b) \llcorner A = 0 \in \Hom(\otimes^{r-1} \lieg , V)$$

The \emph{Killing generators} at $b \in B$ are 
$$ \Kill^\infty(b) = \bigcap_m \Kill^m(b)$$
\end{definition}

Note that $\Kill^m(b)$ is a subspace of $\lieg$ for all $ m \in \NN \cup \{
\infty \}$.  Moreover,
$$ \Kill^m(bp^{-1}) = (\Ad p)(\Kill^m(b))$$ 
Then define
$$ k_m(x) = \dim \Kill^m(b) \qquad k(x) = \dim \Kill^\infty(b)$$
for any $b \in \pi^{-1}(x)$.

Note that for each $m$, the function $k_m(x)$ is lower semicontinuous---that
is, each $x \in B$ has a neighborhood $U$ with $k_m(y) \leq k_m(x)$ for all
$y \in U$.  The same is true for $k(x)$.

The goal is to show that $m$-related points, for $m$ sufficiently large, are actually related by local automorphisms, and that Killing generators of sufficiently high order give rise to local Killing fields.

\begin{definition}
A \emph{local automorphism} between points $b$ and $b'$ of $B$ is a
diffeomorphism $f$ from a neighborhood of $b$ to a neighborhood of $b'$ such
that $f^*\omega = \omega$.  A \emph{local automorphism} between $x$ and $x'$
in $M$ is a diffeomorphism from a neighborhood $U$ of $x$ to a neighborhood
$U'$ of $x'$ inducing an isomorphism of the Cartan geometries $(U,\pi^{-1}(U),\omega)$ and $(U',\pi^{-1}(U'),\omega)$.
\end{definition}

\begin{definition}
A \emph{local Killing field} near $b \in B$ is a vector field $\widetilde{A}$ defined on a
neighborhood of $b$ such that the flow along $\widetilde{A}$, where it is
defined, preserves $\omega$.  A \emph{local Killing field} near $x \in
M$ is a vector field $A$ near $x$ such that the
flow $\varphi^t_A$ along $A$, if it is defined on a neighborhood $U \times (-\epsilon,\epsilon)$ of $(x,0)$, gives an isomorphism of the restricted Cartan
geometry on $U$ with the restricted Cartan geometry on $\varphi^t_A(U)$ for all $t \in (- \epsilon,\epsilon)$.
\end{definition}

Note that a local automorphism between $x,x' \in M$ lifts to a local automorphism from any $b \in \pi^{-1}(x)$ to some $b' \in \pi^{-1}(x')$.  Similarly, a local Killing field near $x \in M$ lifts to a local Killing field near any $b \in \pi^{-1}(x)$.  Local automorphisms and Killing fields on $B$ also descend to $M$; further, the resulting correspondences are bijective, as the next two propositions show.

\begin{proposition}
\label{prop.autloc.inj}
Let $f$ be a nontrivial local automorphism between points $b$ and $b'$ of $B$.  Then $f$ descends to a nontrivial local automorphism $\bar{f}$ between $\pi(b)$ and $\pi(b')$ in $M$.
\end{proposition}

\begin{Pf}
Denote by $P^0$ the identity component of $P$.
In order to ensure that $f$ commutes with $P$, we assume it is defined on a connected neighborhood $U$ of $b$ with $U \cap U P \subseteq U P^0$.  The $P^0$-action is generated by flows along the $\omega$-constant vector fields $X^\ddag$ with $X \in \liep$.  Because $f_*$ preserves all $\omega$-constant vector fields, it commutes with the $P^0$-action.  
Now there is a well-defined extension of $f$ to $U P$ with $f(qp) = f(q) p$ for any $q \in B, p \in P$.  Note that the extended $f$ still preserves $\omega$: if $q \in U$, $p \in P$, then
\begin{eqnarray*}
 \omega_{f(qp)} \circ f_{*qp} & = & \omega_{f(q)p} \circ (R_p)_* \circ f_{*q} \circ (R_p)^{-1}_* \\
& = & (\Ad p^{-1}) \circ \omega_{f(q)} \circ f_{*q} \circ (R_p)^{-1}_* \\
& = & (\Ad p^{-1}) \circ \omega_q \circ (R_p)^{-1}_* \\
& = & \omega_{qp}
\end{eqnarray*}
Now $f$ descends to a diffeomorphism $\bar{f}$ on $\pi(U) \subset M$, and this diffeomorphism is a local automorphism carrying $\pi(b)$ to $\pi(b')$.   

Suppose that $\bar{f}$ were the identity on $\pi(U)$.  Then $f$ would have the form 
$$f(b) = b \cdot (\rho\circ \pi)(b)$$
for $\rho : \pi(U) \rightarrow P$.  Let $N$ be the subgroup of $P$ generated by the image of $\rho$.  We will show $N$ is a normal subgroup of $G$ contained in $P$, contradicting the global assumptions on $G$ and $P$.

On one hand, $f^* \omega = \omega$, while also
$$ (f^* \omega)_b = (\Ad \circ \rho \circ \pi)(b)^{-1} \circ \omega_b + (\rho \circ \pi)_*$$

(see \cite[3.4.12]{sharpe}).
Then for any $Y \in \lieg$ and $x \in \pi(U)$,
$$ Y = ((\Ad \circ \rho)(x))^{-1} Y + (\rho \circ \pi)_* Y$$
So $(\Ad g)(Y) - Y \in \lien$, the Lie algebra of $N$, for all $g \in N$.  Since $G$ is connected, it follows that $hgh^{-1} \in N$ for all $h \in G$, $g \in N$.
\end{Pf}

Similarly, because local Killing fields in $B$ commute with $\omega$-constant vector fields, they commute with the $P^0$-action and descend to $M$.  The local Killing fields near $b \in B$ or $x \in M$ are finite-dimensional vector spaces, and
will be denoted $\Kill^{loc}(b)$ and $\Kill^{loc}(x)$, respectively.  Let
$$ l(x) = \dim \Kill^{loc}(x)$$.

\begin{proposition}
\label{prop.killloc.inj}
For $x = \pi(b)$, 
$$ \Kill^{loc}(b) \cong \Kill^{loc}(x)$$
Moreover, each $x \in M$ has a neighborhood $U_x$ such that $l(y) \geq l(x)$
for all $y \in U_x$.
\end{proposition}

\begin{Pf}
It was observed above that a local Killing field near $x$ lifts to a unique
local Killing field near any $b \in \pi^{-1}(x)$, and it is clear that this
map is linear.  It was also noted above that local Killing fields on $B$
descend to $M$.  This map is linear, and it is injective by an argument
essentially the same as that in the proof of proposition \ref{prop.autloc.inj}
above.  The desired isomorphism follows.

To prove the second statement of the proposition, take a countable nested
sequence of neighborhoods $U_i$ of $x$ with $\cap_i U_i = \{ x \}$.  Let
$\Kill^{loc}_i(x)$ be the subspace of local Killing fields defined on $U_i$.
Because $\Kill^{loc}_i(x) \subseteq \Kill^{loc}_{i+1}(x)$ and $\cup_i
\Kill^{loc}_i(x) = \Kill^{loc}(x)$, these subspaces eventually stabilize to
the finite-dimensional space $\Kill^{loc}(x)$.  Set $U_x = U_i$ once $\Kill^{loc}_i(x) = \Kill^{loc}(x)$.  

For any $y \in U_x$, every $A \in \Kill^{loc}(x)$ determines an element of
$\Kill^{loc}(y)$.  If $A \in \Kill^{loc}(x)$ has trivial germ at $y$, then the
lift $\widetilde{A}$ to $B$ has trivial germ at any $b \in \pi^{-1}(y)$, in
which case it is trivial everywhere
it is defined.  Thus the map $\Kill^{loc}(x) \rightarrow \Kill^{loc}(y)$ is injective for all $y \in U_x$, so $l(y) \geq l(x)$.
\end{Pf}

\begin{proposition}
\label{prop.frobenius.autloc}
Let $(M,B,\omega)$ be real-analytic.  For
any compact $L \subset B$, there exists $m = m(L)$ such that whenever $b,b'
\in L$ are $m$-related, then there is a local automorphism sending $b$ to $b'$.
\end{proposition}

\begin{Pf}
For each $m \geq 1$, denote by $\mathcal{R}^m$ the $C^\omega$ subset of $B
\times B$ consisting of pairs $(b,b')$ with $D^mK(b) = D^mK(b')$. Note that
$\mathcal{R}^{m+1} \subseteq \mathcal{R}^m$.  By the Noetherian property of
analytic sets, there exists $m = m(L)$ such that $\mathcal{R}^k \cap (L \times
L) = \mathcal{R}^m \cap (L \times L)$ for all $k \geq m$.  

Now let $b,b' \in L$ be $m$-related, so they are in fact $\infty$-related.  Define a map $f$ from an exponential neighborhood of $b$ to a neighborhood of $b'$ by
$$ f (\exp_b Y) = \exp_{b'} Y$$
Note that $f(b) = b'$ and $(f^*\omega)_{b} = \omega_{b}$.  For $Y \in \lieg$, denote by $\widetilde{Y}$ the corresponding $\omega$-constant vector feld on $B$.  Now $f$ is a local automorphism if for all $X,Y \in \lieg$ and sufficiently small $t$,
$$ f_*(\widetilde{Y}(\exp(b,tX))) = \widetilde{Y}(\exp(b',tX))$$
This equation is equivalent to 
$$\ddszero \log_{b'} \circ f (\varphi^s_{\widetilde{Y}} \varphi^t_{\widetilde{X}} b) = \ddszero \log_{b'}(\varphi^s_{\widetilde{Y}} \varphi^t_{\widetilde{X}} b') $$

Because $M$ is $C^\omega$, it suffices to show that for all $k \geq 0$,

\begin{eqnarray}
\kddtzero \ddszero \log_{b'} \circ f (\varphi^s_{\widetilde{Y}} \varphi^t_{\widetilde{X}} b) = \kddtzero \ddszero \log_{b'}(\varphi^s_{\widetilde{Y}} \varphi^t_{\widetilde{X}} b')  \label{eqn.derivs}
\end{eqnarray}

By the BCH formula (proposition \ref{prop.bch}), the right-hand side is
$$\kddtzero \ddszero \frac{1}{(k+1)!} \omega_{b'}(a_{k+1}(t\widetilde{X}, s \widetilde{Y}))$$
Each $a_{k+1}(t X, sY)$ is a sum of $(k+1)$-fold brackets of $X$ and $Y$ with coefficients $t^i s^{k+1-i}/c_i$, where $i$ is the multiplicity of $X$, and $c_i$ an integer.  Then 
$$ \kddtzero \ddszero a_{k+1}(tX,sY) = \frac{k!}{c_k} [X, \ldots, X,Y]$$
and the right-hand side of equation (\ref{eqn.derivs}) is
$$ \frac{1}{(k+1) \cdot c_k} \omega_{b'} [ \widetilde{X}, \ldots, \widetilde{X},\widetilde{Y}]$$ 
where $\widetilde{X}$ appears $k$ times in the iterated bracket.

The left-hand side can be written
$$  \kddtzero \ddszero (\log_{b'} \circ f \circ \exp_{b}) \circ \log_b
(\varphi^s_{\widetilde{Y}} \varphi^t_{\widetilde{X}} b) =  \kddtzero \ddszero \log_b (\varphi^s_{\widetilde{Y}} \varphi^t_{\widetilde{X}} b)$$

which, by the BCH formula again, equals
$$ \frac{1}{(k+1) \cdot c_k} \omega_b [\widetilde{X}, \ldots, \widetilde{X},\widetilde{Y}]$$

So it remains to show that these brackets are the same when $b$ and $b'$ are $\infty$-related.  The following lemma completes the proof.
\end{Pf}

\begin{lemma}
\label{lemma.coeffs.curvature}
Let 
$$\Delta_k(b) = [X, \ldots, X,Y] - \omega_b [\widetilde{X}, \ldots, \widetilde{X},\widetilde{Y}]$$
where $X$ occurs $k$ times in each iterated bracket.
Then $\Delta_k$ obeys the recursive formula for all $k \geq 1$
 \begin{eqnarray*}
  \Delta_{k+1}(b) & = & K_b(X,[X,\ldots, X,Y] - \Delta_k(b) ) -  (\widetilde{X}.\Delta_k)(b) + [X, \Delta_k(b)] 
 \end{eqnarray*}

If $b$ and $b'$ are $\infty$-related, then 
$$(\widetilde{X}^r.\Delta_k)(b) = (\widetilde{X}^r.\Delta_k)(b') \qquad \mbox{for all} \ r \geq 0$$ 

If $A$ is a Killing generator at $b$, then
$$ (\widetilde{A}.\widetilde{X}^r.\Delta_k)(b) = 0 \qquad \mbox{for all} \ r \geq 0$$
\end{lemma}

\begin{Pf}
We begin with the recursive formula for $\Delta_k$ when $k=1$:
\begin{eqnarray*}
\Delta_1(b) & = & [X,Y] - \omega_b[\widetilde{X},\widetilde{Y}] \\
& = & K_b(X,Y)
\end{eqnarray*}

For any $r \geq 0$ and $\infty$-related $b$ and $b'$,
\begin{eqnarray*}
(\widetilde{X}^r.\Delta_1)(b) & = & (\widetilde{X}^r.K)_b(X,Y) \\
& = & (D^{r}K_b(X, \ldots, X))(X,Y) \\
& = & (D^rK_{b'}(X, \dots, X))(X,Y) \\
& = & (\widetilde{X}^r.\Delta_1)(b')
\end{eqnarray*}
where $X$ occurs $r$ times in $(X, \ldots, X)$.

Similarly, if $A$ is a Killing generator at $b$, then 
$$(\widetilde{A}.\widetilde{X}^r.\Delta_1)(b) = ((D^{r+1}K_b \llcorner A)(X, \ldots, X))(X,Y)= 0$$

Next suppose the recursive formula for $\Delta_k$ holds up to step $k$.  At the next step,
\begin{eqnarray*}
\Delta_{k+1}(b) & = & [X,[X, \ldots, X,Y]] - \omega_b[\widetilde{X},[\widetilde{X}, \ldots,
\widetilde{X},\widetilde{Y}]]  \\
& = & [X,[X, \ldots, X,Y]] - \omega_b[\widetilde{X},\widetilde{[X, \ldots, X, Y]}] +
\omega_b[\widetilde{X},\omega^{-1} \circ \Delta_k] \\
& = & K_b(X,[X,\ldots, X,Y]) + \omega_b[\widetilde{X},\omega^{-1} \circ \Delta_k] \\
& = & K_b(X,[X,\ldots, X,Y]) - K_b(X,\Delta_k(b)) + (\widetilde{X}. \Delta_k)(b) + [X,\Delta_k(b)] \\
& = & K_b(X,[X, \ldots,X,Y] - \Delta_k(b)) + (\widetilde{X}.\Delta_k)(b) + [X,\Delta_k(b)] 
\end{eqnarray*}
as desired. 



Suppose $(\widetilde{X}^r.\Delta_k)(b) = (\widetilde{X}^r.\Delta_k)(b')$ for
all $r \geq 0$.  Compute
\begin{eqnarray*}
 (\widetilde{X}^r.\Delta_{k+1})(b) & = &  \widetilde{X}^r.(K(X,[X, \ldots, X, Y] - \Delta_k))(b)  \\
& + & (\widetilde{X}^{r+1}.\Delta_k)(b) + [X,(\widetilde{X}^r.\Delta_k)(b)] 
\end{eqnarray*}

Compute inductively
\begin{eqnarray*}
\widetilde{X}^r.(K(X,[X,\ldots,X,Y] - \Delta_k))(b) & = & (\widetilde{X}^r.K)_b(X,[X,\ldots, X,Y] - \Delta_k(b)) \\
& - & \sum_{i=1}^r (\widetilde{X}^{r-i}.K)_b(X, (\widetilde{X}^i.\Delta_k)(b))
\end{eqnarray*}

By the induction hypothesis on $\widetilde{X}^i.\Delta_k$, and because $b$ and $b'$ are $\infty$-related, each term in the above sum is the same at $b$ as at $b'$.  Therefore
$$ \widetilde{X}^r.(K(X,[X,\ldots, X,Y] - \Delta_k))(b) = \widetilde{X}^r.(K(X,[X,\ldots, X,Y] - \Delta_k))(b')$$

and
\begin{eqnarray*}
(\widetilde{X}^r.\Delta_{k+1})(b) & = & \widetilde{X}^r.(K(X,[X, \ldots, X, Y] - \Delta_k))(b') \\
& + & (\widetilde{X}^{r+1}.\Delta_k)(b') +  [X,(\widetilde{X}^r.\Delta_k)(b')] \\
& = & (\widetilde{X}^r.\Delta_{k+1})(b')
\end{eqnarray*}

We leave to the reader the verification that if $(\widetilde{A}.\widetilde{X}^r.\Delta_k)(b) = 0$ for all $r \geq 0$ and $A$ is a Killing generator at $b$, then
$$ (\widetilde{A}.\widetilde{X}^r.\Delta_{k+1})(b) = 0 \qquad \mbox{for all} \ r \geq 0$$
\end{Pf}

Here is the analogue of proposition \ref{prop.frobenius.autloc} relating Killing generators and local Killing fields.

\begin{proposition}
\label{prop.frobenius.killing}
Suppose that $(M,B,\omega)$ is real-analytic.  Then for all $b \in B$, there
exists $m = m(b)$ such that each Killing generator of order $m$ at $b$
determines a unique local Killing field near $b$.
\end{proposition}

\begin{Pf}
The subspaces $\Kill^m(b)$ eventually stabilize, so there is $m = m(b)$ such
that $\Kill^r(b) = \Kill^\infty (b)$ for all $r \geq m$.  Let $A \in
\Kill^\infty(b)$, so $D^rK(b) \llcorner A = 0$ for all $r \geq 1$.

Let $\widetilde{A}(b) = \omega_b^{-1} A$.  Now define $\widetilde{A}$ near $b$ by flowing
along $\omega$-constant vector fields: let
$$\widetilde{A}(\varphi_{\widetilde{Y}}^t b) = \varphi_{\widetilde{Y}*}^t (\widetilde{A}(b))$$
This vector field is
well-defined in an exponential neighborhood of $b$.  Further, for all $Y \in
\lieg$, the bracket $[\widetilde{A},\widetilde{Y}](b) = 0$.  

To show that $[\widetilde{A},\widetilde{Y}] = 0$ in a neighborhood of $b$ for all $Y \in \lieg$, it suffices to show 
$$ (\log_b)_* \left( [\widetilde{A},\widetilde{Y}](\exp(b,tX)) \right) = 0$$
for all $X \in \lieg$ and $t$ sufficiently small.  Because $M$ is
$C^\omega$, it suffices to show
$$ \kddtzero (\log_b)_* \left( [\widetilde{A},\widetilde{Y}](\exp(b,tX)) \right) = 0$$
for all $k \geq 0$.  

As in the proof of \ref{prop.frobenius.autloc} above, this equation follows from the BCH formula and lemma \ref{lemma.coeffs.curvature}.  The reader is invited to refer to the proof of theorem \ref{thm.smooth.regular} and to complete the present proof.
\end{Pf}

\begin{theorem}
\label{thm.can.frobenius}
Let $(M,B,\omega)$ be a compact $C^\omega$ Cartan geometry modeled on $G/P$.
There exists $m \in \NN$ such that any Killing generator at any $b \in B$ of order $m$ gives rise to a unique local Killing field around $b$.  
\end{theorem}

\begin{Pf}
Recall that $\Kill^m(b p^{-1}) = (\Ad p)(\Kill^m(b))$.  Then proposition
\ref{prop.frobenius.killing} above, together with proposition \ref{prop.killloc.inj},
implies that for all $x \in M$, there exists $m(x)$ such that any Killing
generator of order $m(x)$ at any $b \in \pi^{-1}(x)$ determines a local
Killing field near $x$ in $M$.  

Let $U_x$ be the neighborhood given by
proposition \ref{prop.killloc.inj}, on which all local Killing fields near $x$
can be defined.  Shrink $U_x$ if necessary so that $k_{m(x)}(y) \leq k_{m(x)}(x)$ for
all $y \in U_x$.  We wish to show that $m(y) =
m(x)$.  First,
$$ l(x) \leq l(y) \leq k_{m(x)}(y) \leq k_{m(x)}(x)$$
But $l(x) = k_{m(x)}(x)$, so $l(y) = k_{m(x)}(y)$.  A local Killing field $\widetilde{A}$ near $b \in
\pi^{-1}(y)$ is determined by the value
$\omega(\widetilde{A}(b))$, so $\Kill^{loc}(b)$ maps injectively to
$\Kill^m(b)$ for any $m$.  If these spaces have the same dimension for $m = m(x)$, then this map is an
isomorphism---in other words, every Killing generator of order $m(x)$ at any $b
\in \pi^{-1}(y)$ gives rise to a local Killing field near $y$, and $m(y) = m(x)$.

Now take $U_{x_1}, \ldots, U_{x_n}$ a finite subcover of the covering of $M$ by the
neighborhoods $U_x$.  Set $m = \mbox{max}_i \ m(x_i)$.
\end{Pf}

It is well-known that for any $C^\omega$ manifold $B$ equipped with a $C^\omega$ framing, a local Killing field for the framing near any $b_0 \in B$ can be extended uniquely along curves emanating from $b_0$ (see \cite{amores.killing}).  The same is then true in the base of a $C^\omega$ Cartan geometry $M$, because any local Killing field near $x_0 \in M$ has a unique lift to $B$, and local Killing fields in $B$ project to local Killing fields in $M$.  If two local Killing fields of $M$ have the same germ at a point, then they coincide on their common domain of definition.  It follows that if $M$ is simply connected, then extending a local Killing field along curves from some $x_0$ gives rise to a well-defined global Killing field on $M$.  Then we have the following corollary.

\begin{corollary}
\label{cor.can.global}
Let $(M,B,\omega)$ be a compact, simply connected $C^\omega$ Cartan geometry.  There exists $m
\in \NN$
such that for all $b \in B$, every Killing generator at $b$ of order $m$
gives rise to a unique global Killing field on $M$, which in turn gives rise to a $1$-parameter flow of automorphisms of $M$.
\end{corollary}

\section{Stratification theorem in analytic case}


The \emph{$\Kill^{loc}$-relation} is the equivalence relation on
$M$ with $x \sim y$ if $y$ can be reached from $x$ by flowing along a finite sequence
of local Killing fields.  The \emph{$\Kill^{loc}$-orbits} are the equivalence
classes for the $\Kill^{loc}$-relation.  The next result describes the
configuration of these orbits in $M$; it is a version of Gromov's
stratification theorem for compact $C^\omega$ Cartan geometries.

The Rosenlicht stratification theorem says that when an algebraic group $P$ acts algebraically on
a variety $W$, then there exist
$$ U_0 \subset \cdots \subset U_k = W$$
such that $U_i$ is Zariski open and dense in $\cup_{j \geq i} U_j$ and the quotient $U_i
\mapsto U_i / P$ is a submersion onto a smooth algebraic variety (see \cite{rosenlicht.brasil}, \cite[2.2]{gromov.rgs}).

Let $W = \Hom(\otimes^m \lieg, V)$, and define $\Phi : B \rightarrow W$ to be the $P$-equivariant map sending $b$ to the
$\omega$-derivative $D^mK(b)$.  When $(M,B, \omega)$ is algebraic type, the Rosenlicht stratification of $W$ gives
rise to a $\Kill^{loc}$-stratification of $M$.  Recall that a \emph{simple} foliation on a manifold $V$ is one in which the leaves are the fibers of a submersion from $V$ to another manifold $U$.

\begin{theorem}
\label{thm.can.stratification}
Let $(M,B,\omega)$ be a $C^\omega$ Cartan geometry of algebraic type modeled on $G/P$.  Suppose
that $M$ is compact.  Then
there exists a stratification by $\Kill^{loc}$-invariant sets
$$ V_0 \subset \cdots \subset V_k = M$$
such that each $V_i$ is open and dense in $\cup_{j \geq i} V_j$, and the $\Kill^{loc}$-orbits in $V_i$ are leaves of a simple foliation.
\end{theorem}

\begin{Pf}
Let $m$ be given by theorem \ref{thm.can.frobenius}, so that every Killing
generator of order $m$ on $B$ gives rise to a local Killing field on $M$.
Take $V_i = \pi(\Phi^{-1}(U_i))$, where $U_i$ are the pieces of the Rosenlicht
stratification for the $P$-action on the Zariski closure of $\Phi(B)$ in
$W = \Hom(\otimes^m \lieg, V)$.  Then $\cup V_i = M$ and each $V_i$ is open in
$\cup_{j \geq i} V_j$.  Since $\Phi$ is analytic and each $\cup_{j \geq i}
U_j$ is Zariski closed, $\cup_{j > i} V_j$ is an analytic subset of $\cup_{j
  \geq i} V_j$.  Therefore, $V_i$ is also dense in $\cup_{j \geq i} V_j$.

The map $\Phi$ descends to $\bar{\Phi} : M
\rightarrow W / P$.  Each quotient $U_i / P = X_i$ is a smooth variety.  There is the following commutative
diagram.
$$
\begin{array}{ccc}
B & \stackrel{\Phi}{\rightarrow} & W \\
\downarrow &     & \downarrow \\

M & \stackrel{\bar{\Phi}}{\rightarrow} & W/P \\

\cup  &                   & \cup \\
V_i & \rightarrow & X_i
\end{array}
$$

The fibers of the submersion $V_i \rightarrow X_i$ are analytic submanifolds, and the components of the fibers of $\bar{\Phi}$
foliate $V_i$.  Let $X_i'$ be the leaf space of this foliation.  The map $X_i' \rightarrow X_i$ is a local homeomorphism, so $X_i'$ admits the structure of a smooth manifold for which the quotient map $V_i \rightarrow X_i'$ is a submersion.

Now it remains to show that the leaves of these foliations---that is, the
components of the fibers of $\bar{\Phi}$---are $\Kill^{loc}$-orbits.  Let
$\mathcal{F} = \Phi^{-1}(w) \subset B$ for $w \in W$.  Note that $\mathcal{F}
\mapsto \pi(\mathcal{F})$ is a principal bundle, with fiber $P(w)$, the
stabilizer in $P$ of $w$.  For $\bar{w}$ the
projection of $w$ in $W/P$, each component of
$\bar{\Phi}^{-1}(\bar{w})$ in $M$ is the image under $\pi$ of a component of $\mathcal{F}$.

If each component $\mathcal{C}$ of $\mathcal{F}$ is a $\Kill^{loc}$-orbit in
$B$, then each component $\pi(\mathcal{C})$ is a $\Kill^{loc}$-orbit in $M$.
The tangent space $T_b \mathcal{C} = \omega_b^{-1}(\Kill^m(b))$ for all $b \in \mathcal{C}$.
On the other hand, 
$$\omega_b^{-1} (\Kill^m(b)) = \{ X(b) \ : \ X \in \Kill^{loc}(b) \}$$
Thus the $\Kill^{loc}$-orbit of $b$ is contained in $\mathcal{C}$.  A point
$b \in \mathcal{C}$ has a neighborhood $N_b \subset \mathcal{C}$ such that any
$b'
\in N_b$ equals $\varphi_Y^1 b$ for some $Y \in \Kill^{loc}(b)$.  Then given
$a \in \mathcal{C}$, connect $b$ to $a$ by a path and cover this path with
finitely many such neighborhoods to reach $a$ from $b$ by flowing along
finitely many local Killing fields.
\end{Pf}


 

\section{Gromov representation}
\label{section.gromovrep}

Let $(M,B,\omega)$ be a compact $C^\omega$ Cartan geometry of algebraic type
modeled on $G/P$.  The Frobenius theorem gives local Killing fields from
Killing generators of sufficiently high order.  A slight extension of the main
theorem of \cite{bfm.zimemb} gives Killing generators of sufficiently high
order in $\liep$ from big groups $H < \Aut M$.  This latter theorem is a version of Zimmer's embedding theorem---\cite{zimmer.lorentz}, \cite[5.2.A]{gromov.rgs}---in the setting of Cartan geometries.

Combining local Killing fields that arise from the embedding theorem with
certain Killing fields from $H$ gives rise to local Killing fields that
centralize $\lieh$ in theorem \ref{thm.centralizer} (compare
\cite[5.2.A2]{gromov.rgs}, \cite[4.3]{zimmer.gromovrep}).  Local Killing
fields that centralize $\lieh$ lift to the universal cover of $M$ and extend
to global Killing fields.  The fundamental group $\Gamma$ of $M$ preserves
this centralizer $\liec$, and the representation of $\Gamma$ on $\liec$ is
related to the adjoint representation of $H$ in theorem \ref{thm.gromovrep}, a
version of Gromov's representation theorem \cite[6.2.D1]{gromov.rgs}.  In our
centralizer theorem, the group $H < \Aut M$ is not assumed to preserve a
finite volume.  In neither the centralizer nor the representation theorem is it assumed simple; see \cite{nz.smooth.factors} for some related statements on existence of Gromov representations for simple $H$ without a finite invariant measure, in the setting of Gromov's rigid geometric structures.

\subsection{Embedding theorem}

If $H$ is a Lie subgroup of $\Aut M$, then the Lie algebra $\lieh$ can be viewed as an algebra of global Killing fields on $B$.  If $b \in B$ and $X$ is a nontrivial Killing field on $B$, the evaluation $X(b) \neq 0$.  There are therefore for each $b \in B$ linear injections $\iota_b : \lieh \rightarrow \lieg$ defined by 
$$\iota_b(X) = \omega_b(X)$$

The embedding theorem relates the adjoint representation of $H$ on $\lieh$ with the representation of a certain subgroup of $P$ on $\iota_b(\lieh)$.  The key ingredient in the proof of the embedding theorem is the Borel density theorem.  It essentially says that a finite measure on a variety that is invariant by an algebraic action of a group $S$ is supported on $S$-fixed points.  One must take care, however, that $S$ has no nontrivial compact quotients.

\begin{definition}
Let $H$ be a Lie group.  A Lie subgroup $S < H$ is \emph{discompact} if the Zariski closure $\mbox{Zar}(\mbox{Ad}_\lieh S)$ has no nontrivial compact algebraic quotients. 
\end{definition}


The following statement is a consequence of the Borel density theorem and appears in \cite[3.2]{bfm.zimemb}.

\begin{theorem}(see \cite[2.6]{dani.borel.density} and \cite[3.11]{shalom.discompact})
\label{thm.bd}
Let $\psi : S \rightarrow \Aut W$ for $S$ a locally compact group and $W$ an algebraic variety, and assume that $\mbox{Zar}(\psi(S))$ has no nontrivial compact algebraic quotients.  Suppose $S$ acts continuously on a topological space $M$ preserving a finite Borel measure $\mu$.  Assume $\phi : M \rightarrow W$ is an $S$-equivariant measurable map.  Then $\phi(x)$ is fixed by $\mbox{Zar}(\psi(S))$ for $\mu$-almost-every $x \in M$.
\end{theorem}

Now we can state the embedding theorem that will be needed.

\begin{theorem}
\label{thm.embedding}
Let $(M,B,\omega)$ be a Cartan geometry of algebraic type modeled on $G/P$.
Let $H < \Aut M$ be a Lie group and $S < H$ a discompact subgroup preserving a
probability measure $\mu$ on $M$.  Denote by $\bar{S}$ the Zariski closure of
$\mbox{Ad}_\lieh S$.  For any $m \geq 0$, there exists $\Lambda \subset B$
with $\mu( M \setminus \pi(\Lambda)) = 0$, such that to every $b \in \Lambda$ corresponds an algebraic subgroup $\check{S}_b < \mbox{Ad}_\lieg P$ with
\begin{enumerate}
\item{$\check{S}_b (\iota_b(\lieh)) = \iota_b(\lieh)$}
\item{the representation of $\check{S}_b$ on $\iota_b(\lieh)$ is equivalent to $\bar{S}$ on $\lieh$}
\item{$\check{S}_b$ fixes $D^rK(b)$ for all $0 \leq r \leq m$}
\end{enumerate}
\end{theorem}

The proof is the same as in \cite{bfm.zimemb}, except that we apply theorem \ref{thm.bd} to strata in the $P$-quotient of the variety 
$$ \widetilde{W} = \Mon(\lieh,\lieg) \times V \times \cdots \times \Hom(\otimes^m \lieg,V)$$
the target of the $P$-equivariant map
$$\widetilde{\phi}(b) = (\iota_b, K(b), \ldots, D^mK(b))$$ 

The action of $\bar{S}$ on $W$ is by 
$$ g(\rho, \varphi_0, \ldots, \varphi_m) = (\rho \circ g^{-1}, \varphi_0,
\ldots, \varphi_m)$$
The action of $p \in P$ on $\Mon(\lieh,\lieg)$ is by post-composition with
$\Ad p$.  Then $p$ acts on the first factor of $W$ by this action, and on the
remaining factors by the
actions defined in section \ref{section.frobenius} above: for $\varphi \in
\Hom(\otimes^m \lieg,V)$, 
$$ p.\varphi = p \circ \varphi \circ \Ad^m p^{-1}$$
Note $\tilde{\phi}$ is $P$-equivariant.    

Let $\phi : M \rightarrow W/P$ be the map induced by $\widetilde{\phi}$; it is
$S$-equivariant.  By theorem \ref{thm.bd}, for $\mu$-almost-every $x$, the point $\phi(x)$ is fixed by $\bar{S}$.  Let $x$ be a such a point, and let $b \in \pi^{-1}(x)$.  Then define
$$ \check{S}_b = \{ p \in \AdgP \ : \ p.\widetilde{\phi}(b) = g.\widetilde{\phi}(b) \ \mbox{for some} \ g \in \bar{S} \}$$

Then $\check{S}_b$ satisfies the conditions (1)-(3) of theorem \ref{thm.embedding}.

\subsection{Centralizer theorem}

The group $\check{S}_b$ gives rise, via the Frobenius theorem, to elements of the stabilizer of $\pi(b)$, which in turn give local Killing fields commuting with $\lieh$.  Denote by $\widetilde{M}$ the universal cover of $M$ and by $q$ the covering map.  Denote by $\mathfrak{c}$ the Lie algebra of global Killing fields on $\widetilde{M}$ commuting with the algebra $\lieh$ of Killing fields lifted from the $H$-action on $M$.  Let $\lies$ be the Lie algebra of $S$.  Given a point $y$ of a manifold $N$ and an algebra $\lieu$ of vector fields, denote by $\lieu(y)$ the subspace of $T_yN$ consisting of values at $y$ of elements of $\lieu$.

\begin{theorem}
\label{thm.centralizer}
Let $(M,B,\omega)$ be a compact $C^\omega$ Cartan geometry of algebraic type.  Let $H < \Aut M$ be a Lie group and $S < H$ a discompact subgroup preserving a probability measure $\mu$ on $M$.  Then for $\mu$-almost-every $x \in M$, for every $\tilde{x} \in q^{-1}(x)$, the subspace $\lies(\tilde{x}) \subset \liec(\tilde{x})$.
\end{theorem}

\begin{Pf}
The ideas of the proof are the same as Zimmer's \cite{zimmer.gromovrep}.  Let
$m$ be given by theorem \ref{thm.can.frobenius}, so that any Killing
generator of order $m$ at any $b \in B$ gives rise to a unique local Killing
field.  Let $x$ belong to the full-measure set $\Lambda$ as in the embedding theorem \ref{thm.embedding}, and let $b
\in \pi^{-1}(x)$.  Denote by $\check{\lies}$ the Lie algebra of $\check{S}_b$,
and by $\rho_b$ the Lie algebra homomorphism of $\check{\lies}$ onto
$\bar{\lies}$, the Lie algebra of $\bar{S}$.
For any $X \in \check{\lies}$ and $0 \leq r \leq m$, 
$$ D^r K(b) \llcorner X = (X^\ddag.D^{r-1}K)(b) = - X.(D^{r-1}K(b)) = 0$$
by $P$-equivariance of $D^rK$.  Therefore $\check{\lies} \subset \Kill^m(b)$.  Now the Frobenius theorem guarantees, for each $X \in \check{\lies}$, a local Killing field $X^*$ near $b$ with $\omega_b (X^*) = X$.  Because $X^*(b)$ is tangent to the fiber over $x$, the local Killing field near $x$ induced by $X^*$ fixes $x$.

Let $Y \in \lieh$, viewed as a Killing field on $B$.  Compute
\begin{eqnarray*}
\mbox{d} \omega (X^*,Y) & = & X^*. \omega(Y) - (L_Y\omega)(X^*) \\
& = & X^*.\omega(Y) \\
& = & (L_{X^*} \omega)(Y) + \omega[X^*,Y] \\
& = & \omega[X^*,Y]
\end{eqnarray*}

On the other hand, since $\omega_b (X^*) \in \liep$, the curvature $\Omega_b(X^*,Y) = 0$, so 
\begin{eqnarray*}
 \omega_b[X^*,Y] & = &  \mbox{d} \omega_b(X^*,Y)  \\
& = & [\omega_b(X^*) ,\omega_b(Y)]  = [X,\iota_b (Y)] \\
& = & \iota_b ((\rho_b X)(Y)) = \omega_b ((\rho_b X)(Y))
\end{eqnarray*}

Both $[X^*,Y]$ and $(\rho_b X)(Y)$ are local Killing fields.  They are
determined by their values at any point of $B$, so they must be equal.  We conclude that for all $Y \in \lieh$,
$$ [X^*,Y] = (\rho_b X)(Y)$$

Now, given $X \in \lies$ and $x \in M$ satisfying the conclusion of the
embedding theorem, choose any $b \in \pi^{-1}(x)$ and let $Y^*$ be the local
Killing field on $M$ fixing $x$ with $(\rho_b \circ \omega_b)( Y^*) = \mbox{ad}_{\lieh} X$.  Then define $X^c = X - Y^*$.  It is a local Killing field near $x$ satisfying
\begin{itemize}
\item{ $X^c(x) = X(x) - Y^*(x) = X(x)$}
\item{ for all $W \in \lieh$,
$$[X^c,W] = [X - Y^*,W] = [X,W] - ((\rho_b \circ \omega_b)(Y^*))(W) = 0$$ } 
\end{itemize}

Now $X^c$ lifts to a local Killing field near any $\tilde{x} \in
\widetilde{M}$.  Because $\widetilde{M}$ is real-analytic and simply
connected, there is a unique global extension of $X^c$ to $\widetilde{M}$, which will also be denoted $X^c$.  Now $X^c \in \liec$, and $X^c(\tilde{x}) = X(\tilde{x})$.  Such an $X^c$ exists for any $X \in \lies$, so the theorem is proved.
\end{Pf}


\subsection{Gromov representation}

We first review Zimmer's notion of the algebraic hull of a measurable cocycle.
Two references on this subject are \cite{zimmer.mackey} and \cite{feres.dsss}.

\begin{definition}
Let $S$ be a locally compact group acting on a topological space $M$ preserving an ergodic probability measure $\mu$.  Let $L$ be a topological group. An $L$-valued \emph{measurable cocycle} for the $S$-action on $M$ is a measurable map $\alpha : S \times M \rightarrow L$ satisfying
$$ \alpha(gh,x) = \alpha(g,hx) \alpha(h,x)$$
for all $g,h \in S$ and almost-every $x \in M$.
\end{definition}

\begin{definition}
Let $S$ be a locally compact group acting by automorphisms of a $V$-vector
bundle $E$ over a topological space $M$.  Suppose that $S$ preserves an
ergodic probability measure $\mu$ on $M$.  A \emph{measurable trivialization}
of $E$ is a measurable map $t : E  \rightarrow M \times V$ of the form $t(x,v)
= (x, t_xv)$, where $t_x$ is a linear isomorphism $E_x \rightarrow V$ for
almost-every $x$.
\end{definition}

A measurable trivialization $t$ gives rise to a $\GL(V)$-valued measurable
cocycle $\alpha_t$ where
$$ t(g(x,v)) = (gx,\alpha(g,x) (t_xv))$$

\begin{definition}
Let $S, M, \mu, V, E$ be as in the previous definition.  The \emph{algebraic hull} of the $S$-action is the minimal algebraic subgroup $L < \GL(V)$ for which there exists a measurable trivialization $t$ of $E$ with $\alpha_t(S \times M) \subseteq L$.
\end{definition}

The algebraic hull is well defined up to conjugacy in $\GL(V)$; this is a consequence of the Borel density theorem. See \cite{zimmer.mackey}.

We will need the following fundamental facts about the algebraic hull.  A \emph{virtual epimorphism} of algebraic groups is a homomorphism $\sigma : L_1 \rightarrow L_2$ for which $\sigma(L_1)$ is a Zariski dense subgroup of $L_2$ of finite index.

\begin{proposition}
\label{prop.alghull}
Let $S, M, \mu, V, E$ be as above.  
\begin{enumerate}
\item{Let $\widetilde{M}$ be the universal cover of $M$, $\Gamma \cong
    \pi_1(M)$, and $\widetilde{S}$ the connected group of lifts of $S$ to $\widetilde{M}$.  Let $\rho : \Gamma \rightarrow \GL(V)$ be a representation and let $E = \widetilde{M} \times_\rho V$.  Then $\widetilde{S}$ acts by automorphisms of $E$, and the algebraic hull is contained in $\Zar(\rho(\Gamma))$.}

\item{Let $E_0$ be an $S$-invariant subbundle of $E$.  There is a virtual epimorphism from the algebraic hull of $S$ on $E$ to the algebraic hull of $S$ on $E_0$.}

\item{Let $E_0$ be as above, and let $E' = E/E_0$.  There is a virtual epimorphism from the algebraic hull of $S$ on $E$ to the algebraic hull of $S$ on $E'$.}

\item{Suppose there is a trivialization $t$ of $E$ in which $\alpha_t(g,x) = \rho(g)$ for $\rho : S \rightarrow \GL(V)$ a homomorphism.  Then the algebraic hull of the $S$-action is $\Zar(\rho(S))$.}
\end{enumerate}
\end{proposition}

\begin{Pf}
For (1), note that $\widetilde{S}$ commutes with $\Gamma$, so the $\widetilde{S}$-action on $\widetilde{M} \times V$ by $g(\tilde{x},v) = (g \tilde{x}, v)$ commutes with the $\Gamma$-action on the product.  Then the $\widetilde{S}$-action on $\widetilde{M} \times V$ descends to $E$.  The rest is proposition 3.4 of \cite{zimmer.gromovrep}; it is a straightforward exercise with measurable cocycles.

Items (2), (3), and (4) are straightforward; they appear as propositions 3.3 and 3.5 of \cite{zimmer.gromovrep}.
\end{Pf}

Let $S < \Aut M$ be as above.  Denote by $\lies_x$ the Lie algebra of the
stabilizer in $S$ of $x \in M$.  Suppose that $\lies_x$ is an ideal $\lies_0
\lhd \lies$.  Then denote by $J(S,x) = \Zar(\bar{\mbox{Ad }}S)$, where $\bar{\mbox{Ad}}$ is the representation of $S$ on $\lies / \lies_0$ obtained as a quotient of the adjoint representation.  

\begin{theorem}
\label{thm.gromovrep}
Let $(M,B,\omega)$ be a compact $C^\omega$ Cartan geometry of algebraic type.  Let $S < \Aut M$ be discompact, and suppose that $S$ preserves a probability measure $\mu$ on $M$.  Then for $\mu$-almost-every $x \in M$, $\lies_x \lhd \lies$, and there is a representation $\rho$ of $\pi_1(M) \cong \Gamma$ for which $\Zar(\rho(\Gamma))$ contains a subgroup with a virtual epimorphism to $J(S,x)$.
\end{theorem}

\begin{Pf}
By decomposing $\mu$ into ergodic components if necessary, we may assume that $\mu$ is ergodic.

We first present the standard argument due to Zimmer that almost every stabilizer is an ideal.  Define a map 
\begin{eqnarray*}
\psi & : & M \rightarrow \Gr \lies = \bigcup_{k=0}^{\dim \lies} \mbox{Gr}^k \lies \\
     &   & x \mapsto \lies_x
\end{eqnarray*}

The group $S$ acts on $W = \Gr \lies$ via $\Ad S$, and $\Zar(\Ad S)$ has no compact algebraic quotients by the discompactness assumption.  The map $\psi$ is $S$-equivariant.  The Borel density theorem \ref{thm.bd} thus applies, and for $\mu$-almost-every $x$, the stabilizer $\lies_x$ is $\Ad S$-fixed---in other words, it is an ideal $\lies_0$.

Now suppose $\lies_x = \lies_0 \lhd \lies$ and in addition that $x$ satisfies
the conclusion of the centralizer theorem \ref{thm.centralizer}, so
$\lies(\tilde{x}) \subset \liec(\tilde{x})$ for every $\tilde{x} \in
q^{-1}(x)$.  Because the Killing fields of $\lies$ on $\widetilde{M}$ are
lifted from $M$, they commute with $\Gamma$.  Therefore the centralizer
$\liec$ is normalized by $\Gamma$.  Let $\rho$ be the representation of
$\Gamma$ on $\liec$.  By proposition \ref{prop.alghull} (1), the algebraic
hull of $\widetilde{S}$ on $E = \widetilde{M} \times_\rho \liec$ is contained
in $\Zar(\rho(\Gamma))$.  Note that in fact the $\widetilde{S}$-action on $E$
factors through $S$, because any element of $\widetilde{S} \cap \Gamma =
\ker(\widetilde{S} \rightarrow S)$ centralizes $\liec$.

Denote by $T\mathcal{O}$ the tangent bundle to $S$-orbits in $M$
$$T \mathcal{O} = \{ (x, Y(x)) \ : \ x \in M, \ Y \in \lies \}$$
There is an obvious measurable trivialization $t : T \mathcal{O} \rightarrow M
\times \lies / \lies_0$ in which the cocycle for the $S$-action is
$\alpha(g,x) = \bar{\mbox{Ad }} g$.  Then by proposition \ref{prop.alghull} (4), the algebraic hull of $S$ on $T\mathcal{O}$ equals $J(S,x)$.

The evaluation map $\epsilon : \widetilde{M} \times \liec \rightarrow T\widetilde{M}$ with $\epsilon(\tilde{x},Y) = Y(\tilde{x})$ descends to an $S$-equivariant map $\bar{\epsilon} : E \rightarrow TM$.  The kernel $E_0$ is an $S$-invariant subset of $E$, in which each fiber $(E_0)_x$ is a vector subspace of $E_x$.  The dimension of $(E_0)_x$ is $S$-invariant, so we may consider $E_0$ a subbundle of $E$.   The algebraic hull of $S$ on $E$ virtually surjects onto the algebraic hull of $S$ on $E' = E / E_0$ by proposition \ref{prop.alghull} (3).

The map $\epsilon$ factors through an isomorphism almost-everywhere from $E' = E/E_0$ to an $S$-invariant subbundle $\bar{\epsilon}(E)$ of $TM$, so the algebraic hulls on these two are isomorphic.   But $\bar{\epsilon}(E)$ also contains the $S$-invariant subbundle $T \mathcal{O}$, so the algebraic hull of $S$ on $\bar{\epsilon}(E)$ virtually surjects onto the algebraic hull of $S$ on $T\mathcal{O}$ by proposition \ref{prop.alghull} (2).

We conclude that the algebraic hull of $S$ on $E$, which is contained in $\Zar(\rho(\Gamma))$, virtually surjects onto $J(S,x)$, as desired.
\end{Pf}

\begin{corollary}
Let $S < \Aut M$ be semisimple with no compact local factors.  Suppose that $S$ preserves a finite volume form on $M$. Then there is a representation $\rho$ of $\pi_1(M) \cong \Gamma$ for which $\Zar(\rho(\Gamma))$ contains a subgroup with a virtual epimorphism to $\Zar(\Ad S)$.
\end{corollary}

\begin{Pf}
Let $\mu$ be the finite measure determined by the $S$-invariant volume form on
$M$.  There are only finitely-many nontrivial ideals of $\lies$.  For each nonzero ideal $\lies_x$, the fixed set has empty interior (see
\cite[7.1]{bfm.zimemb}).  The $S$-action is thus locally
free---that is, $\lies_x = 0$---almost everywhere.  Then $J(S,x) = \Zar(\Ad S)$.  Since also $S$ is
discompact, the corollary follows from theorem \ref{thm.gromovrep}.
\end{Pf}

\section{Frobenius and open-dense results in smooth case}

The analytic Frobenius theorem says that a Killing generator at any $x \in M$
gives rise to a local Killing field.  In this section we show that Killing
generators of smooth Cartan geometries still give rise to local Killing fields on an open dense subset of $M$, consisting of the \emph{regular} points.  Recall that $k(x)$, for $x \in M$ is the dimension of $\Kill^\infty(b)$ for any $b \in \pi^{-1}(x)$.

\begin{definition}
Let $(M,B,\omega)$ be a $C^\infty$ Cartan geometry.  The \emph{regular} points of $M$ are those $x \in M$ for which $k(x)$ is locally constant.
\end{definition}

Because $k(x)$ is lower semicontinuous, the regular points are an open, dense subset of $M$.

\begin{proposition}
\label{prop.extn.reg.curves}
 Suppose that, for $X \in
\lieg$, the curve $\gamma(t) = \exp(b,tX)$ consists of regular points for all
$t \in [-1,1]$.   Then there exists $m$ such that
$\Kill^m(\gamma(t)) = \Kill^\infty(\gamma(t))$ for all $t \in [-1,1]$. Moreover, for any $b \in
B$ and $A \in \Kill^\infty(b)$,
$$\omega_{\gamma(t)} (\varphi_{\widetilde{X}*}^tA) \in \Kill^\infty(\gamma(t))$$
for all $t \in [- 1,1]$. 
\end{proposition}

\begin{Pf}
Let $m(b)$ be such that $k_r(b) = k(b)$ for all $r \geq m(b)$.  For all $t$ sufficiently small, 
$$ k(\gamma(t)) \leq k_{m(b)}(\gamma(t)) \leq k_{m(b)}(b) = k(b)$$
The regularity assumption means $k(\gamma(t)) = k(b)$, so $k_{m(b)}(\gamma(t)) = k(\gamma(t))$ for all $t$ sufficiently small.  Now repeating the argument along the compact curve $\gamma$ shows that $k_{m(b)}(\gamma(t)) = k(\gamma(t))$ for all $t \in [-1,1]$, and $\Kill^r(\gamma(t)) = \Kill^\infty(\gamma(t))$ for all $r \geq m(b)$.

Let as above $V = \wedge^2(\lieg/\liep)^* \otimes \liep$, and for $r \in \NN$, let 
$$\mathcal{W}^r = \bigoplus_{i=0}^{r} \Hom(\otimes^i \lieg, V)$$  

where we set $\otimes^0 \lieg = \RR$.  For $X \in \lieg$ and $(K_0, \ldots,
K_r) \in \mathcal{W}^r$, write
$$ (K_0, \ldots, K_r) \llcorner X = (K_1 \llcorner X, \ldots, K_r \llcorner X) \in \mathcal{W}^{r-1}$$

Now denote as usual by $K$ the curvature function $B \rightarrow V$.  For $b
\in B$, let 
$$ \mathcal{D}^rK(b) = (K(b), \ldots, D^rK(b)) \in \mathcal{W}^r$$
and 
\begin{eqnarray*}
\mathcal{C}^r_b & : & \lieg \rightarrow \mathcal{W}^{r-1} \\
&  & X \mapsto \mathcal{D}^rK(b) \llcorner X
\end{eqnarray*}

The kernel of $\mathcal{C}^r_b$ is $\Kill^r(b)$.  By the discussion in the previous
paragraph, $\ker \mathcal{C}^{m(b)}_{\gamma(t)} = \ker
\mathcal{C}^{m(b)+1}_{\gamma(t)}$ for all $t$.  Therefore, the functionals on
$\lieg$ appearing in
the decomposition of $\mathcal{C}^{m(b)+1}_{\gamma(t)}$ in terms of a basis of $\mathcal{W}^{m(b)}$ are linear combinations of the functionals appearing in any
decomposition of $\mathcal{C}^{m(b)}_{\gamma(t)}$ in terms of any basis of
$\mathcal{W}^{m(b)-1}$.


Denote $\omega_{\gamma(t)}(\varphi_{\widetilde{X}*}^t A) = A(t)$ and by $\widetilde{A}$ the corresponding vector field along $\gamma$.  Then, for each $1 \leq r \leq m(b)$,
\begin{eqnarray*}
 \frac{\mbox{d}}{\mbox{d}t} \left( D^rK(\gamma(t))\llcorner A(t) \right) & = & \left( \widetilde{X}.\widetilde{A}.D^{r-1}K \right)(\gamma(t)) \\
& = & \left( \widetilde{A}.\widetilde{X}.D^{r-1}K \right)(\gamma(t)) \\
& = & \left( D^{r+1}K(\gamma(t) ) \llcorner A(t) \right)\llcorner X
\end{eqnarray*}

using that $[\widetilde{X},\widetilde{A}] = 0$.
There results a system of ODEs
$$ \frac{\mbox{d}}{\mbox{d}t}  \mathcal{C}^r_{\gamma(t)} (A(t)) =
\mathcal{C}^{r+1}_{\gamma(t)} (A(t)) \llcorner X$$
as $r$ ranges from $1$ to $m(b)$.
At $t=0$, all $\mathcal{C}^r_{\gamma(0)} (A(0)) = 0$.  Then 
$$\mathcal{C}^r_{\gamma(t)} (A(t)) = \mathcal{D}^rK(\gamma(t))\llcorner A(t) \equiv 0$$
 is the unique solution for all $1 \leq r \leq m(b)+1$, and $A(t) \in \Kill^{m(b)}(\gamma(t)) = \Kill^\infty(\gamma(t))$ for all $t$. 
\end{Pf}

\begin{theorem}
\label{thm.smooth.regular}
Let $(M,B,\omega)$ be a $C^\infty$ Cartan geometry and let $U \subseteq M$ be
the set of regular points.  For each
component $U_0 \subseteq U$, there exists $m = m(U_0)$ such that every Killing
generator of order $m$ at any $b \in \pi^{-1}(U_0)$ gives rise to a unique
local Killing field near $\pi(b)$.
\end{theorem}

\begin{Pf}
Let $b \in U_0$, and let $m$ be such that $\Kill^m(b) = \Kill^\infty(b)$.
Then by proposition \ref{prop.extn.reg.curves}, for all $b' \in U_0$, there is
also $\Kill^m(b') = \Kill^\infty(b')$.  So it suffices to show that any
Killing generator at a point lying over the regular set determines a local
Killing field.

Let $A \in \Kill^\infty(b)$ for $b \in \pi^{-1}(U)$.  As in the proof of proposition \ref{prop.frobenius.killing}, we define a vector field $\widetilde{A}$ in an exponential neighborhood of $b$ by $\widetilde{A}(\exp(b,tX)) = \varphi_{\widetilde{X}*}^tA$.  By proposition \ref{prop.extn.reg.curves}, $\omega(\widetilde{A})$ is a Killing generator everywhere it is defined.

To show that $\widetilde{A}$ descends to a local Killing field near $\pi(b)$, it suffices to show it is a local Killing field near $b$.  Then we must show that for any $Y,X \in \lieg$ and sufficiently small $T$, the bracket 
$$ [\widetilde{A},\widetilde{Y}](\exp(b,TX)) = 0$$

We will show that, in the chart $\log_b$, this field satisfies the ODE
$$  \frac{\mbox{d}}{\mbox{d}t} \log_{b*} \left( [\widetilde{A},\widetilde{Y}](\exp(b,tX)) \right) = 0 $$
Because the initial value at $t=0$ is zero, this will imply vanishing for all $t$.

Let $b(T) = \exp(b,TX)$ and $\Psi_T = (\log_b \circ \exp_{b(T)})_*$.  Then
\begin{eqnarray*} 
\left. \frac{\mbox{d}}{\mbox{d}t} \right|_T  \log_{b*} \left( [\widetilde{A},\widetilde{Y}](\exp(b,tX)) \right) & = &
\ddtzero  \log_{b*} \left( [\widetilde{A},\widetilde{Y}](\exp(b,(T+t)X)) \right) \\
& = & \ddtzero \left[ (\Psi_T \circ \log_{b(T)*}) \left( [\widetilde{A},\widetilde{Y}](\exp(b(T),tX)) \right) \right]
\end{eqnarray*}

So it suffices to show that for each $T$,
$$ \ddtzero \log_{b(T)*} \left( [\widetilde{A},\widetilde{Y}](\exp(b(T),tX)) \right) = 0$$
Now
 \begin{eqnarray*}
\log_{b(T)*} \left( [\widetilde{A},\widetilde{Y}](\exp(b(T),tX)) \right) & = & \ddszero
\log_{b(T)*} \left( \varphi_{\widetilde{Y}*}^s
  (\widetilde{A}(\varphi_{\widetilde{X}}^t b(T)))
-  \widetilde{A}(\varphi_{\widetilde{Y}}^s \varphi_{\widetilde{X}}^t b(T)) \right) \\
& = & \ddszero \log_{b(T)*} \left( \varphi_{\widetilde{Y}*}^s
  \varphi_{\widetilde{X}*}^t (\widetilde{A}(b(T))) - \varphi_{\widetilde{Z}(t,s)*}^1
  (\widetilde{A}(b(T))) 
  \right) 
\end{eqnarray*}

where 
$$Z(t,s) = (\log_{b(T)} \circ \varphi^s_{\widetilde{Y}} \circ
\varphi^t_{\widetilde{X}}) (b(T)) = \zeta_{b(T)}(tX,sY)$$ 
as in the BCH formula.  Write $\widetilde{A}_T = \widetilde{A}(b(T))$.  Now
\begin{eqnarray}
& & \log_{b(T)*} \left( \varphi_{\widetilde{Y}*}^s \varphi_{\widetilde{X}*}^t
  (\widetilde{A}_T) - \varphi_{\widetilde{Z}(t,s)*}^1
  (\widetilde{A}_T)   
\right)  = \\
& & \label{eqn.log.bracket}
\dduzero \left[ ( \log_{b(T)} \circ \varphi_{\widetilde{Y}}^s \circ \varphi_{\widetilde{X}}^t)
(\varphi_{\widetilde{A}}^u b(T)) \right] - (\log_{b(T)} \circ \varphi_{\widetilde{Z}(t,s)}^1)_*
(\widetilde{A}_T) 
\end{eqnarray}

Let $c(u) = \varphi^u_{\widetilde{A}}b(T)$.  The first term of (\ref{eqn.log.bracket}) can be written 
\begin{eqnarray}
& & \dduzero \left[ \left( (\log_{b(T)} \circ \exp_{c(u)}) \circ
(\log_{c(u)} \circ \varphi_{\widetilde{Y}}^s \circ
\varphi_{\widetilde{X}}^t ) \right) (c(u)) \right] \\
& = & \label{eqn.log.flow.comp}
\left[ \dduzero (\log_{b(T)} \circ \exp_{c(u)})
\right](Z_0(t,s)) + \dduzero Z_u(t,s)
\end{eqnarray}

where $Z_0(t,s) = Z(t,s)$, and 
$$ Z_u(t,s) = (\log_{c(u)} \circ
\varphi_{\widetilde{Y}}^s \circ
\varphi_{\widetilde{X}}^t)(c(u)) =
\zeta_{c(u)}(tX,sY)$$

Now the first term of (\ref{eqn.log.flow.comp}) is
\begin{eqnarray*}
\left[ \dduzero (\log_{b(T)} \circ \exp_{c(u)})
\right](Z_0(t,s)) 
& = & \dduzero (\log_{b(T)} \circ
\varphi_{\widetilde{Z}_0(t,s)}^1 \circ \varphi_{\widetilde{A}}^u)(b(T)) \\
& = & (\log_{b(T)} \circ \varphi_{\widetilde{Z}_0(t,s)}^1)_*(\widetilde{A}_T)
\end{eqnarray*} 

Thus the first term of (\ref{eqn.log.flow.comp}) cancels with the second term of
(\ref{eqn.log.bracket}), and we are left to show
$$ \ddtzero \ddszero \dduzero Z_u(t,s)  = \ddtzero \ddszero \dduzero \zeta_{c(u)}(tX,sY)= 0$$

We have
\begin{eqnarray*}
\ddtzero \ddszero \dduzero \zeta_{c(u)}(tX,sY) & = & 
\dduzero \left( \ddtzero \ddszero \zeta_{\varphi_{\widetilde{A}}^u b(T)}
  (tX,sY) \right) \\
& = & \dduzero \left(\frac{1}{2} \omega_{\varphi_{\widetilde{A}}^u b(T)}[\widetilde{X},\widetilde{Y}] \right) \\
& = & \frac{1}{2} \dduzero \left( [X,Y] - K_{\varphi_{\widetilde{A}}^u b(T)}(X,Y) \right) \\
& = & - \frac{1}{2} (\widetilde{A}.K)_{b(T)}(X,Y) \\
& = & 0
\end{eqnarray*}
because $\widetilde{A}(b(T))$ is a Killing generator.
\end{Pf}

\begin{theorem}
\label{thm.smooth.opendense}
Let $(M,B,\omega)$ be a $C^\infty$ Cartan geometry of algebraic type.  Suppose that $M$ contains
a dense $\Kill^{loc}$-orbit.  Then $M$ contains an open, dense, locally
homogeneous subset.
\end{theorem}

\begin{Pf}
Let $\mathcal{O} \subset M$ be a dense
$\Kill^{loc}$-orbit.  Because the regular set $U$ is open and $\Kill^{loc}$-invariant,
it contains $\mathcal{O}$.  Because $\mathcal{O}$ is connected, $U$ has only one component.
Let $m$ be such that for all $b \in
\pi^{-1}(U)$, any Killing generator of order $m$ at $b$ gives rise to a local
Killing field near $\pi(b)$ (such $m$ exists by \ref{thm.smooth.regular}).  

The map $\Phi : B \rightarrow \Hom(\otimes^m \lieg, V)$ gives rise to a
stratification as in \ref{thm.can.stratification}
$$ V_1 \subset \cdots \subset V_k = M$$
such that $\bar{\Phi}$ is a smooth map of each $V_i$ onto a smooth variety.
Because $V_1$ is open and $\Kill^{loc}$-invariant, it contains $\mathcal{O}$.
Therefore, $V_1 \cap U$ is open and dense. 
The same argument as for
\ref{thm.can.stratification} shows that components of fibers of $\bar{\Phi}$ in $V_1
\cap U$ are $\Kill^{loc}$-orbits, and they are closed in $V_1$.  Then
$$\mathcal{O} = \bar{\mathcal{O}} \cap V_1 \cap U = V_1 \cap U $$
so $\mathcal{O}$ is an open, dense, locally homogeneous subset of $M$.
\end{Pf}

\begin{question}
This question is asked in \cite{dag.rgs} section 7.3: Can the conclusion of theorem
\ref{thm.smooth.opendense} above be strengthened to say that $M$ is
locally homogeneous?
\end{question}

The forthcoming corollary gives a positive answer in a very special case.  For
$(M,B,\omega)$ a Cartan geometry modeled on $G/P$, the tangent bundle $TM$ can
be identified with $B \times_P (\lieg /\liep)$ (see \cite[4.5.1]{sharpe}).
The Cartan geometry will be called \emph{unimodular} when the representation
of $P$ on $\lieg / \liep$ has image in $\SL(\lieg / \liep)$.  In this case,
there is a volume form on $(M,B,\omega)$ preserved by $\Aut M$.

\begin{corollary} (see \cite[1.8]{bfm.zimemb})
Let $(M,B,\omega)$ be a compact, simply connected, unimodular, $C^\omega$ Cartan geometry of
algebraic type.
Let $H < \Aut M$ be a connected Lie subgroup.  If $H$ has a dense orbit in $M$, then $M$ is homogeneous: there exists
$H' < \Aut M$ acting transitively.  
\end{corollary}

\begin{Pf}
If $H$ has a dense orbit in $M$, then there is a dense $\Kill^{loc}$-orbit in
$M$.  By theorem \ref{thm.smooth.opendense}, there is an open dense
$\Kill^{loc}$-orbit $U \subseteq M$.  But all local Killing fields on $M$
extend to global ones because $M$ is $C^\omega$ and simply connected (see
\cite{amores.killing}), and they are complete because $M$ is compact.  Then
the volume-preserving automorphism group of $M$ has an open orbit.  The
conclusion then follows from theorem 1.7 of \cite{bfm.zimemb}.
\end{Pf}

\bibliographystyle{ieeetr.bst}
\bibliography{karinsrefs}













\end{document}